\documentclass[12pt]{amsart}

\usepackage{graphicx}
\usepackage{amssymb}

\theoremstyle{plain}
\newtheorem{theorem}{Theorem}
\newtheorem{proposition}{Proposition}[section]
\newtheorem{lemma}[proposition]{Lemma}
\newtheorem{corollary}[proposition]{Corollary}

\theoremstyle{definition}
\newtheorem{definition}[proposition]{Definition}

\theoremstyle{remark}

\newcommand{\eps}[2]{{\hspace{-3pt}\begin{array}{c}%
  \raisebox{-2.5pt}{\includegraphics[width=#1]{#2.eps}}%
\end{array}\hspace{-3pt}}}

\def\qed{{ \hfill $\eps{8mm}{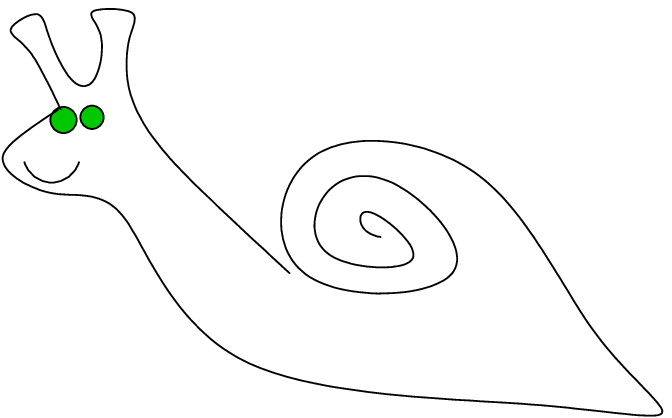}$ }}

\begin{document}

\title[]{on chern-simons theory with an inhomogeneous gauge group and BF
theory knot invariants}

\author{Gad Naot}

\address{University of Toronto, Department of Mathematics, Canada}

\email{gnaot@math.toronto.edu}

\begin{abstract}
We study the Chern-Simons topological quantum field theory with an
inhomogeneous gauge group, a non-semi-simple group obtained from a
semi-simple one by taking its semi-direct product with its Lie
algebra. We find that the standard knot observable (i.e. trace of
the holonomy along the knot) essentially vanishes, and yet, the
non-semi-simplicity of the gauge group allows us to consider a class
of un-orthodox observables which breaks gauge invariance at one
point and leads to a non-trivial theory on long knots in
$\mathbb{R}^3$.
\\
We have two main morals :

1. In the non-semi-simple case there is more to observe in
Chern-Simons theory! There might be other interesting non
semi-simple gauge groups to study in this context beyond our
example.

2.  In the case of an inhomogeneous gauge group, we find that
Chern-Simons theory with the un-orthodox observable is actually the
same as 3D BF theory with the Cattaneo-Cotta-Ramusino-Martellini
knot observable. This leads to a simplification of their results and
enables us to generalize and solve a problem they posed regarding
the relation between BF theory and the Alexander-Conway polynomial.
We prove that the most general knot invariant coming from pure BF
topological quantum field theory is in the algebra generated by the
coefficients of the Alexander-Conway polynomial.

\end{abstract}

\maketitle

\newpage

\section{\bf{Introduction}}
We would like to address the question of the most general knot
invariant coming from BF topological quantum field theory (TQFT). In
the mid-90's, Cattaneo et al. ~\cite{ca1,ca2,ca3,ca4} showed that
while BF theory with cosmological constant produces the same
invariants of knots as the Chern-Simons (CS) theory, the BF theory
with no cosmological constant (pure BF theory) and SU(2) gauge group
produces invariants that lie in the algebra generated by the
coefficients of the Alexander-Conway polynomial. The BF TQFT is
completely equivalent to CS theory, however while the equivalence
with non-zero cosmological constant maintains the semi-simplicity
property of the gauge group, the equivalence when the cosmological
constant is set to zero shifts us to a CS theory with a
non-semi-simple gauge group. Thus, starting with a gauge group G in
the pure BF theory we end up on the CS side with a semi-direct
product of the group with its Lie algebra $G\ltimes \mathcal{G}$,
also known as the inhomogeneous group of G (denoted IG). Eventually,
the question we want to address is as follows:

\begin{center}
\emph{What is the most general knot invariant (or knot observable)
in a CS theory with an inhomogeneous gauge group ?}
\end{center}

\vspace{3mm}The natural thing to begin with is the standard gauge
invariant observable: the trace of the holonomy of the gauge
connection along the knot, in some chosen representation. As we
shall see in chapter 2, this observable gives us no information
regarding the knot (except for framing information that can be
normalized to zero anyway). This fact seems to pose a contradiction
since in the equivalent pure BF theory one can extract some
non-trivial information regarding the knot (at least in the case of
SU(2) gauge). Therefore, there must be a procedure that will give us
non-trivial information regarding the knot in this setting.
\\

In this paper we will introduce the following procedure (section
3.1): We take an observable which is not gauge invariant at one
single point on the knot (holonomy along the knot without the
trace). By doing that, we get a CS theory on $S^3$ with broken gauge
symmetry at one point on the knot. Then, we declare this point to be
infinity by taking the point out of our space (i.e. puncturing $S^3$
at that point). Since the gauge transformations are taken to vanish
at infinity, we get a completely invariant CS theory on
$\mathbb{R}^3$ with a long knot embedded in it. Since knot theory in
$S^3$ and long knot theory in $\mathbb{R}^3$ are ``isomorphic''
theories we lose no information as far as knot theory is concerned
and we get ``legal'' CS theory in which we can consider a new and
wider class of observables.\\

Using perturbation theory, we will show that our observable gives
non-trivial information about the knot (section 3.2). Continuing
with perturbation theory in our setting (CS with IG gauge group and
the new observable) we will be able to prove that the most general
knot invariant coming from such a construction is in the algebra
generated by the coefficients of the Alexander-Conway polynomial
(section 3.3). Returning to BF theory (chapter 4) we will create an
observable that reproduces Cattaneo et al.'s result for SU(2) and
generalizes it to every metrizable gauge group and representation
thereof used in the construction of the pure BF theory.
\\

Our construction detects, in some sense, the difference between the
invariant subspace of the universal enveloping algebra of the
non-semi-simple gauge group and the co-invariant quotient of the
universal enveloping algebra. This enables us to extract non-trivial
information regarding the knot (a difference which does not exist in
the semi-simple case). We discuss this issue in chapter 5.

\section{\bf{Perturbation theory for the standard observable with an inhomogeneous gauge group}}
In this section we will show that the perturbation theory of
Chern-Simons (CS) theory with an inhomogeneous gauge group and the
standard knot observable is almost trivial. We start (2.1) with some
reminders about perturbation theory (in CS context). Unfortunately,
due to the size of the topic this reminder is not meant to teach the
theory. For recent detailed pedagogic reviews see ~\cite{sa,pol}.
This will be followed by introduction of the inhomogeneous gauge
group and our notations for it (2.2). We continue with giving more
details about perturbation theory using this gauge group (2.3).
Section 2.3 will describe the consequences of using such a gauge
group on the knot invariants coming from the standard observable.
Finally, in section 2.4, we will show that the standard construction
in this setting is trivial, leading to (almost) trivial knot
invariants.

\subsection{Some reminders about perturbation theory and Chern-Simons theory}
We recall the setting for knot invariants in the framework of CS
theory (in this chapter a knot is an embedding of $S^1$ into $S^3$).
We take $S^3$ to be the 3 dimensional manifold on which the CS
theory is defined. We choose any gauge group G whose lie algebra is
metrizable (e.g. semi-simple gauge group). Denote the metric on the
Lie algebra by $<,>$ and let A be a G-connection. We consider the
following action functional:
\[
\mbox{Chern-Simons action}\hspace{2mm} CS(A)=k \cdot2\pi \int_{S^3}
<A \wedge dA + \frac{2}{3} A \wedge A \wedge A>
\]
where $k$ is a coupling constant which satisfies a quantization
condition (e.g. for semi-simple gauge groups it must be an integer).
In what follows, however, we will work perturbatively with formal
power series in $k^{-1}$.\\

Now, recall the standard knot observable - $Tr_Rhol(A)$ where
$hol(A)$ stands for the holonomy of the connection A along the knot
(the path-ordered exponent of the integral of A along the knot) and
the trace is taken in some chosen representation $R$ of the
gauge group.\\

One obtains a knot invariant by taking the path integral over all
connections with the observable plugged into the integrand :
\[
Z(knot)=\int \mathcal{DA} ~e^{iCS(A)}Tr_Rhol(A)
\]
This integral is usually referred to as the expectation value of the
standard observable and it is viewed as a function from knots to the
base field (the complex numbers from now on). One is usually
interested in $\frac{Z(knot)}{Z(\emptyset)}$, where $Z(\emptyset)$
is the path
integral with no observable plugged in (no knot).\\

Calculating the integral using perturbation theory we get the
following sum:
\[
\sum_{D \in D_{CS}} {k^{-degD} \cdot \zeta^{CS}(D) \cdot
W_{\mathcal{G}}(D)}
\]
where $D_{CS}$ is defined to be the set of all trivalent connected
graphs based on $S^1$ (Feynman diagrams based on a Wilson loop),
$\zeta^{CS}$ is the integration of the corresponding propagators
over the diagram (that is, the integration over the appropriate
configuration space) and $W_{\mathcal{G}}$ is the Lie-algebraic part
of the expectation value which holds all the information coming from
the Lie group G and its representation $R$ (that is, the weight
system, reviewed quickly below).
\\

In the CS perturbation theory calculation, every internal edge
contributes $1/k$ and every internal vertex contributes $k$ to the
total power of $k$ of the diagram (recall that an internal
edge/vertex is an edge/vertex which is not on the Wilson loop). The
total power of $1/k$ is called the \emph{degree} of the diagram. A
convenient way of counting the degree is labeling each internal
vertex with -1 and each edge with +1 and then summing all the labels
to get the degree of the diagram. This number is denoted
$degD$.\\

We want to factorize the perturbation theory calculation in order to
get a better hold on the above knot invariant. As a first step, we
can drop the weight system and replace it with the diagram D on
which it is calculated. We get the following invariant:
\[
\sum_{D \in D_{CS}} {k^{-degD} \cdot \zeta^{CS}(D) \cdot D}
\]
 This invariant is viewed
as taking values in $\mathcal{A}(S^1)$, where we define
$\mathcal{A}(S^1)$ as follows :

\begin{definition}
$\mathcal{A}(S^1)$ is the algebra of all connected trivalent graphs
based on a circle (i.e. $D_{CS}$) quotient by the IHX relation
$\eps{20mm}{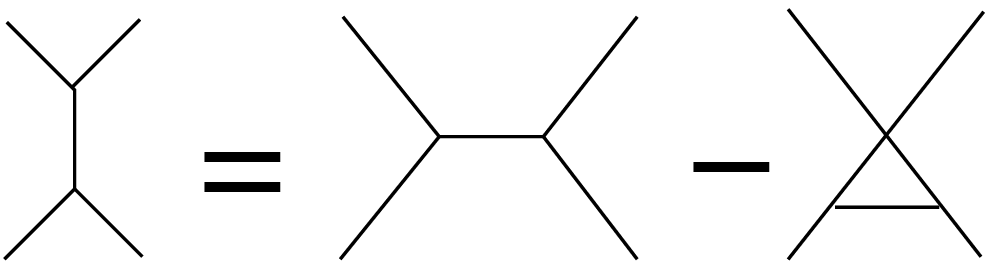}$, the STU relation $\eps{30mm}{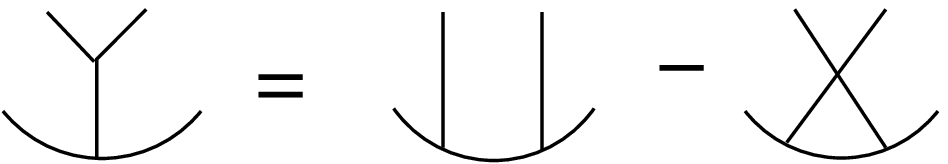}$ and the
anti-symmetry relation $\eps{15mm}{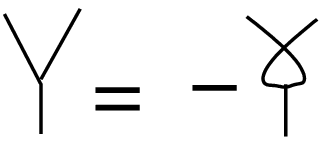}$.
\end{definition}

When we apply the weight system $W_{\mathcal{G}}$ to this invariant
we get the first sum. Next we will apply the weight
system by parts, in order to get a further factorization.\\

Weight systems are described and explored in \cite{ba3}. We recall
the main construction of a weight system coming from a Lie algebra,
which includes three parts: labeling, contracting  and the trace
part. In the labeling part, after choosing a basis for the Lie
algebra, one labels each internal edge of a given diagram (i.e.
edges that are not contained in the base circle) with a different
index on each end of the edge. In the second part, one writes the
structure constants tensor for each internal vertex and the metric
tensor for each internal edge (using the indices on these
vertices/edges and the chosen basis) and proceeds to contract
matching indices (using the metric to raise and lower indices). At
that point, one is left with an invariant tensor in the enveloping
algebra of $\mathcal{G}$ whose indices are the ones on the vertices
on the base loop. In the final part, the trace part, one represents
this tensor using a chosen representation of $\mathcal{G}$ and takes
the trace. This concludes the construction of the weight system.
\\

Taking a diagram (more precisely a class representative but we will
not make these distinctions throughout the paper) in
$\mathcal{A}(S^1)$ and applying to it the labeling and contracting
parts of the weight system gives us a tensor in the universal
enveloping algebra. This tensor is invariant under the adjoint
action (all the tensors used in building it are invariant under the
action and contractions are $\mathcal{G}$-maps) and therefore it is
in the invariant subspace of the universal enveloping algebra
$U(\mathcal{G})^\mathcal{G}$. Although we did not apply the trace
part of the weight system yet, we still have a trace of the trace
that appears in the standard observable we started with, which comes
about diagrammatically in the fact that we are looking at diagrams
based on $S^1$. The ability to move legs around cyclically (which is
the trace property) forces us to quotient the resulting tensor in
$U(\mathcal{G})^\mathcal{G}$ into the co-invariant quotient of the
algebra, thus getting an invariant in
$(U(\mathcal{G})^\mathcal{G})_\mathcal{G}$. Finally, one can apply
the trace part to get elements of $\mathbb{C}$.

\vspace{3mm} To summarize, when we use perturbation theory for the
standard observable, we get a knot invariant in the form of a series
(in $1/k$) with coefficients in $\mathbb{C}$, which factorize in the
following way:

\begin{equation} \{knots\} \rightarrow \mathcal{A}(S^1)[[k^{-1}]] \rightarrow
(U(\mathcal{G})^\mathcal{G})_\mathcal{G}[[k^{-1}]] \rightarrow
\mathbb{C}[[k^{-1}]]
\end{equation}

\vspace{3mm} Although there is much more to say about Chern-Simons
theory and the use of perturbation theory in this context, we will
stop at this point. We recommend \cite{ba1,ba2,ax1,alv1,alt1,gua1}
for various approaches (physically and/or mathematically inclined)
in addition to the reviews cited above. The fact that we actually
get a knot invariant in every step of the factorization above is an
example of an important issue which we completely ignored here.
\\

\subsection{The inhomogeneous gauge group}
We now present a specific type of gauge group that we wish to focus
on in the context of perturbation theory of CS theory.\\

Start with a semi-simple gauge group G and take a semi-direct
product of it with its Lie algebra (the semi-direct action is the
adjoint action). Look at the Lie algebra of this semi-direct product
and denote it $L_0$. As a vector space, $L_0$ is a double copy of
the original algebra, $\mathcal{G} \bigoplus
\overline{\mathcal{G}}$, where we use upper bar notation to
distinguish the second copy. For any $X\in\mathcal{G}$ let $X
\mapsto \overline{X}$ be the identity map between the two copies.
Thus $\overline{X}$ is the element $X$ in the second copy of
$\mathcal{G}$. Let $X_i$ be a basis for $\mathcal{G}$,
$\overline{X_i}$ the corresponding basis for
$\overline{\mathcal{G}}$ and [.,.] the original bracket structure on
$\mathcal{G}$. We have the following bracket structure on $L_0$:
\[
[X_i,X_j]_{L_0} = [X_i,X_j] \hspace{6mm} [X_i,\overline{X_j}]_{L_0}
= \overline{[X_i,X_j]} \hspace{6mm}
[\overline{X_i},\overline{X_j}]_{L_0}=0
\]
We see now that the second copy is the abelianization of
the original algebra.\\

Letting the metric on $\mathcal{G}$ be $<,>$ (the one coming from
the trace of the adjoint representation, say) we take the following
(invariant) metric on $L_0$:
\[
<X_i,X_j>_{L_0} = 0 = <\overline{X_i}, \overline{X_j}>_{L_0}
\hspace{6mm} <X_i,\overline{X_j}>_{L_0}=<X_i,X_j>
\]
\\

\subsection{The structure of the perturbation theory expansion with the inhomogeneous gauge group}
Let us look at the perturbation theory expansion using the
inhomogeneous gauge group $L_0$. As we will show, the factorization
(1) factors even further as the weight system using this type of Lie
algebra is more refined. We will define the space of directed
``legal'' diagrams and construct a factorization
of the weight system $W_{L_0}$ through a map into it.\\

As known from standard perturbation theory of quantum gauge
theories, the algebraic contributions of vertices and edges in a
Feynman diagram are determined by the structure constants and the
metric on the Lie algebra. These contributions are encoded in the
weight system, as described earlier.

Given a diagram, choose a basis $X_i$ to $\mathcal{G}$, take the
corresponding basis $\overline{X_i}$ to $\overline{\mathcal{G}}$
(obtaining a basis to the entire algebra) and apply the weight
system using the following steps:
\\

1(pre). Sum over all the ways of labeling the ends of the edges of
the diagram with different indices $i$ and $\bar{i}$ representing
the basis for $\mathcal{G}$ and $\overline{\mathcal{G}}$
respectively. Therefore, instead of labeling with one type of
indices running over the entire basis of the algebra, label with two
types of indices, each running over one summand in the direct sum
that forms the
algebra.\\

2(pre). Write the structure constants tensor and the metric tensor
(as before) with indices according to the labeling, and contract the
matching indices (thus in each summand of step 1(pre) only a part
of the algebra basis is contracted).\\

Since the summations are finite, one can see that this way of
applying the labeling and contraction parts of the weight system
(i.e. in two summation steps having two types of algebra indices
instead of one) is equivalent to the usual way of applying the
weight system.
\\

Furthermore, we know that certain combinations of indices vanish in
the second step. By looking at the metric we observe that the only
edges which are non zero are the ones that have a non-bar index on
one end and a bar-index at the other end. By looking at the bracket
structure we observe that the only combination of indices on a
vertex that will not be zero at step 2(pre) is a combination of two
non-bar indices and one bar index.
This allows us to revise step 1(pre) above into:\\

1(revised). Sum over all ways of labeling the ends of the edges of
the diagram with indices $i$ and $\bar{i}$, representing the basis
for $\mathcal{G}$ and $\overline{\mathcal{G}}$ respectively, in such
a way that each edge has one bar index and one non-bar index at its
ends, and each vertex has two non-bar indices and one bar index at
its legs.\\

Let us introduce a notation. Wherever we have an edge after step
1(revised) we direct the edge (by putting an arrow head) from the
non-bar index to the bar index and drop the bar on top of the bar
index. Thus the two steps for applying the weight system of $L_0$ on
a diagram can be finalized into the following form :
\\

1. Sum over all ways of directing the edges of the diagram (internal
edges only) in such a way that each vertex has two outgoing legs and
one ingoing leg (we call it a ``legal'' directing and it looks like
$\eps{7mm}{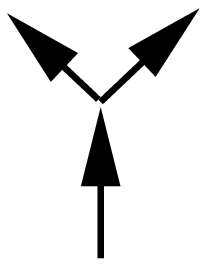}$).

2. Put a different index on each end of all the directed edges,
write down the appropriate tensors and contract indices according to
the arrow convention above (i.e. the arrow's direction indicates if
the index is considered bar or non-bar, and we use the appropriate
basis accordingly).
\\

Step 1 is a map from the space of all Feynman diagrams $D^{CS}$ to
the space of all directed Feynman diagrams
$\overrightarrow{D^{CS}}$, defined to be all trivalent connected
graphs based on a circle with arrows on the internal edges (the
edges which are not a part of the base circle) such that each vertex
is ``legal'' in the above sense. This map takes a diagram to the sum
of all possible ``legal'' directing of it. In order to extend this
map to $\mathcal{A}(S^1)$ we need to treat the STU/IHX/AS relations.
For that purpose we define the directed STU, directed IHX and
directed AS relations which we will use in taking a quotient of
$\overrightarrow{D^{CS}}$. These relations reflect the algebraic
structure of $L_0$ and can be read directly from the bracket
structure.

\begin{definition}
The Directed STU relations in $\overrightarrow{D^{CS}}$ are the
following relations :
\begin{center}
$\begin{array}{llll} \eps{7mm}{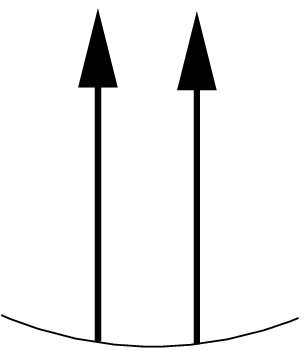}-\eps{7mm}{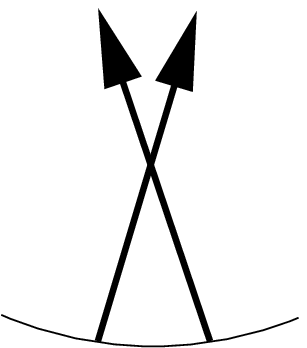}=\eps{7mm}{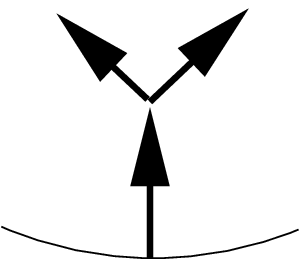} &
& &
\eps{7mm}{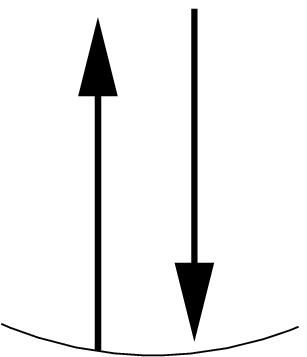}-\eps{7mm}{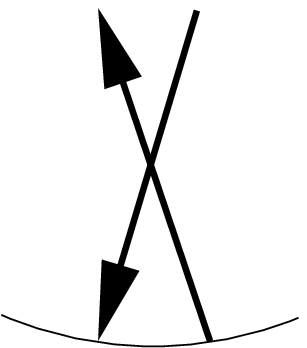}=\eps{7mm}{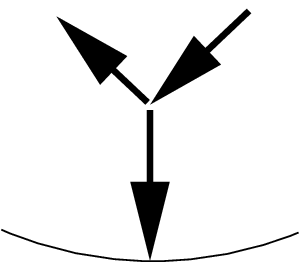} \\
\eps{7mm}{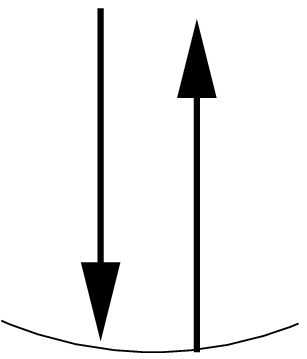}-\eps{7mm}{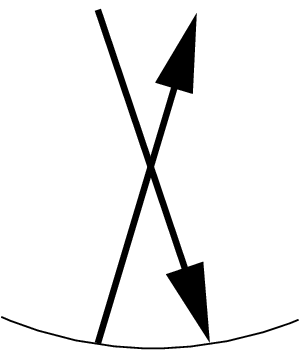}=\eps{7mm}{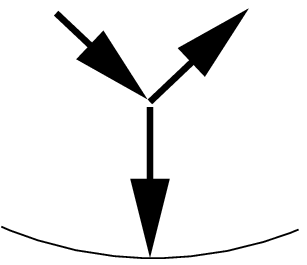} & & &
\eps{7mm}{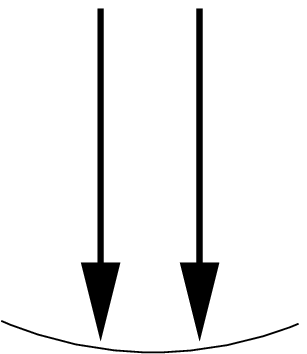}-\eps{7mm}{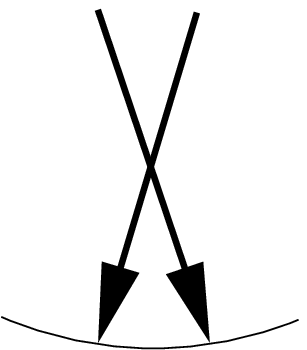}=\eps{7mm}{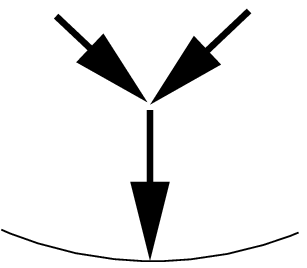}=0
\end{array}$
\end{center}
Where the bottom non-directed arc is a part of the base circle.\\

The Directed IHX relations are actually a consequence of the
directed STU relations in $\overrightarrow{D^{CS}}$, just as in
~\cite{ba3}. We draw one example here and the rest are just
rotations thereof :

\begin{center}
$\eps{5mm}{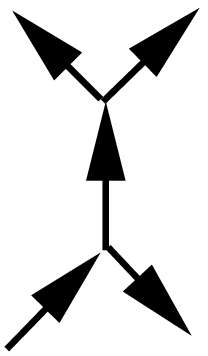} = \eps{7mm}{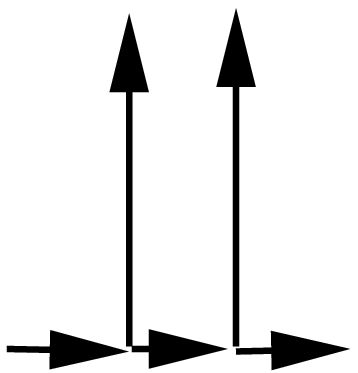} - \eps{7mm}{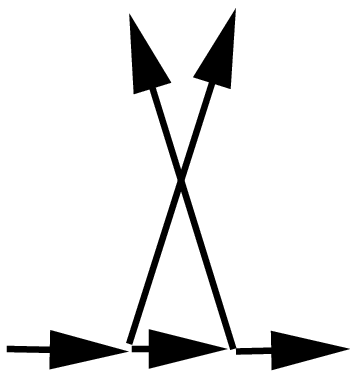}$
\end{center}

The Directed AS relation is defined just as the usual AS relation :

\begin{center}
$\eps{7mm}{as1} + \eps{7mm}{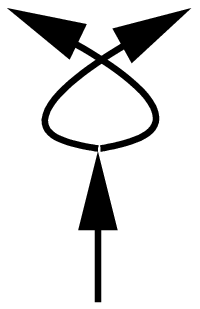} = 0$
\end{center}
\end{definition}

\begin{definition}
The space of all ``legally'' directed diagrams modulo the directed
STU/IHX and AS relations, denoted
$\overrightarrow{\mathcal{A}(S^1)}$, is the space of all trivalent
connected graphs based on a circle (the Wilson loop) with arrows on
the internal edges such that each vertex looks like $\eps{7mm}{as1}$
(i.e. the space $\overrightarrow{D_{CS}}$ defined above), quotient
by the directed STU, directed IHX and directed AS relations defined
above.
\end{definition}

Armed with these definitions we can see that step 1 of the
factorization above is well defined as a map $\mathcal{A}(S^1)
\longrightarrow \overrightarrow{\mathcal{A}(S^1)}$. Composed with
step 2 above (labeling and contracting of indices) and with taking
the trace of the resulting tensor in some chosen representation, we
get the following factorization of the perturbation theory expansion
of CS theory with an inhomogeneous gauge group and the standard
observable :

\begin{equation}
\{knots\} \longrightarrow \mathcal{A}(S^1)[[k^{-1}]] \longrightarrow
\overrightarrow{\mathcal{A}(S^1)}[[k^{-1}]] \longrightarrow
(U(L_0)^{L_0})_{L_0}[[k^{-1}]] \longrightarrow \mathbb{C}[[k^{-1}]]
\end{equation}

\vspace{3mm} Note that the directed STU relations (which reflect the
structure of $L_0$) tell us that two adjacent legs which touch the
base circle and have arrows pointed towards it, are commutative (as
the bar-indices are commutative). On the other hand, if they
originate from one vertex they are anti-commutative according to the
AS relation. Thus one gets an important relation in
$\overrightarrow{\mathcal{A}(S^1)}$:
\[
\eps{10mm}{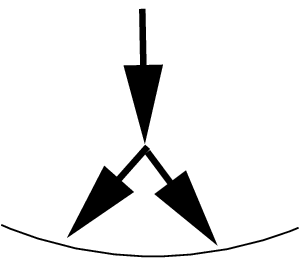} = 0
\]

This relation will be used in the arguments presented in the next
section.

\subsection{But $\vec{\mathcal{A}}(S^1)$ is (almost) empty!}
We will now prove that there are only few non-zero diagrams in
$\overrightarrow{\mathcal{A}(S^1)}$ and the ones that are not zero
encode a very specific type of information regarding the knot ---
the framing. This information can actually be normalized to zero
showing that the entire weight system $W_{L_0}$ is trivial. We will
prove that the only primitive element which is non-zero in
$\overrightarrow{\mathcal{A}(S^1)}$ is the directed single chord
diagram $\eps{4mm}{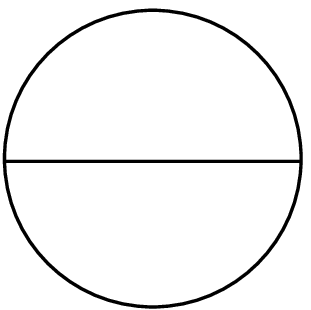}$ (drawn without the arrow here) and thus
the only degree $n$ diagram which is non-zero in
$\overrightarrow{\mathcal{A}(S^1)}$ is $\eps{4mm}{theta}^n$. We will
do that by first showing that the only primitive elements that are
possibly non-zero are the wheel diagrams (e.g.
$\eps{4mm}{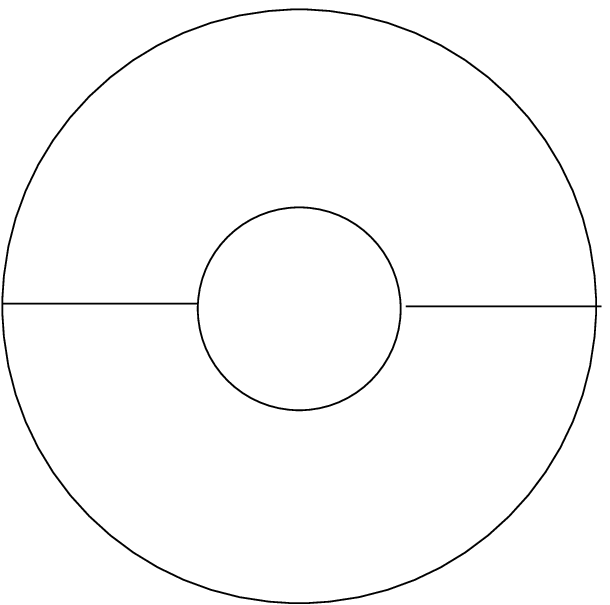},\eps{4mm}{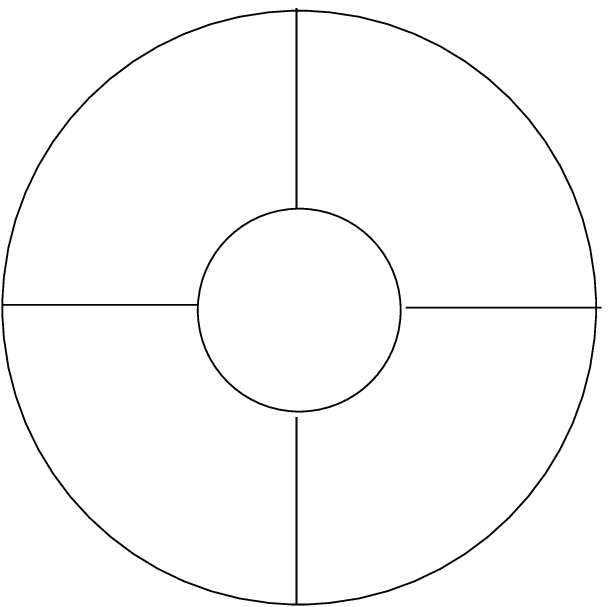}$) and $\eps{4mm}{theta}$. We
will then present an argument as to why the wheel diagrams are
actually zero. We remind the reader that whenever we mention a
diagram we mean a (representative of) diagram class.

\begin{lemma}
Every degree n diagram in $\overrightarrow{\mathcal{A}(S^1)}$ has at
least n external vertices (vertices that are on the base loop).
\end{lemma}

\begin{proof}
As already mentioned, to every diagram we attach a power of $1/k$
which is called the degree. A convenient way of counting the degree
is labeling each internal vertex with -1 and each edge with +1 and
then summing up all the labels to get the degree of the diagram.
Given a degree $n$ diagram in $\overrightarrow{\mathcal{A}(S^1)}$,
label it with +1 and -1 according to the above. ``Push'' the labels
on the edges in the direction of the arrows. Whenever a label hits a
vertex it stops and the labels at the vertex are added. An example
of this procedure for a degree 3 diagram is given by:
\[
\eps{40mm}{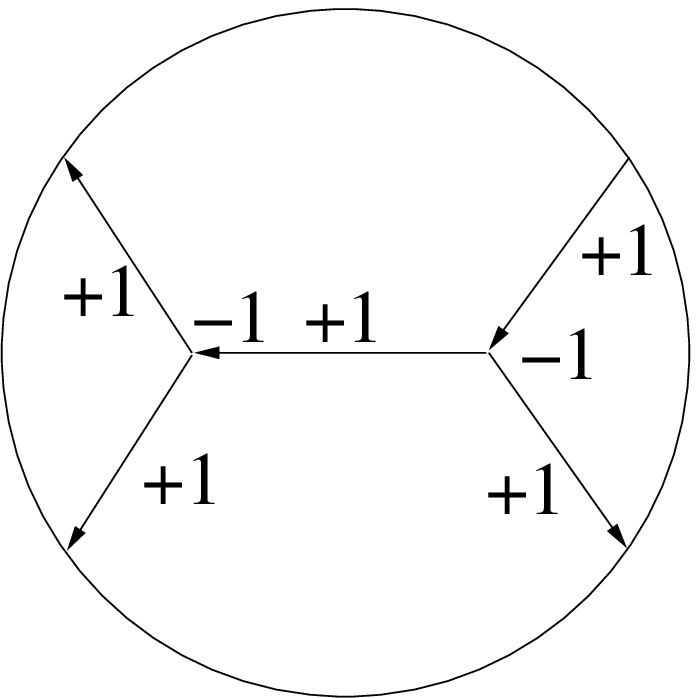} \xrightarrow{push} \eps{40mm}{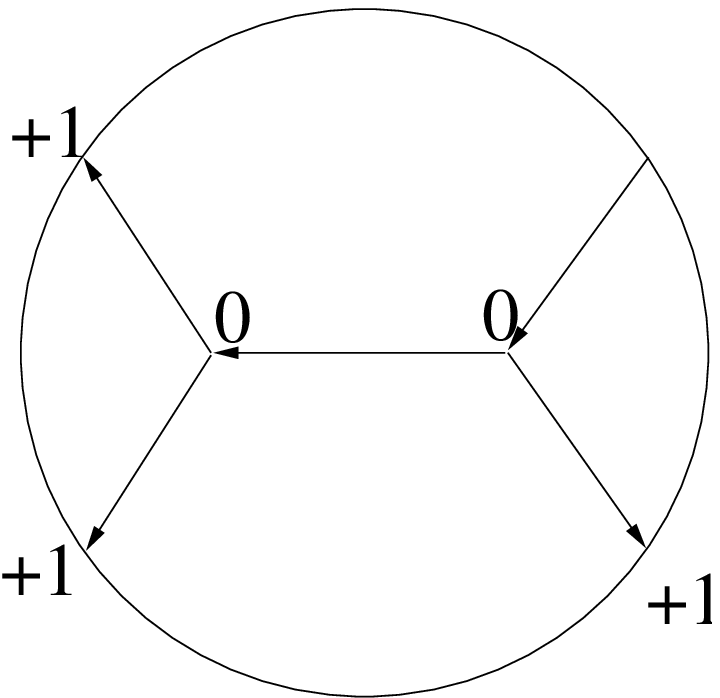}
\]

\vspace{3mm} According to the ``legal'' directing rules, each
internal vertex will have one label +1 hitting the vertex label -1
thus summing to a zero label. All other labels are pushed to the
external vertices (on the Wilson loop). Since the total sum of
labels is conserved and equal to the degree of the diagram, what we
have just proven is that a degree $n$ diagram in
$\overrightarrow{\mathcal{A}(S^1)}$ has at least $n$ vertices on the
Wilson loop (in other words at least $n$ external
legs).\\

\qed
\end{proof}

\vspace{3mm} This lemma has immediate consequences in terms of the
primitive elements that generate the algebra
$\overrightarrow{\mathcal{A}(S^1)}$. The primitive elements of the
non-directed algebra $\mathcal{A}(S^1)$ are these diagrams which
remain connected when the Wilson loop is removed from them. Let
$\mathcal{P}_n$ be the space of primitive diagrams in
$\mathcal{A}(S^1)$ of degree $n$. Filter the primitive spaces (of
different degrees) according to the number of external legs and
write $\mathcal{P}_{n,d}$ for the space of primitive diagrams in
$\mathcal{A}(S^1)$ of degree $n$ and at most $d$ external legs.

\begin{lemma} \cite{chm1,chm2}
If $n$ is even, then $\mathcal{P}_n = \mathcal{P}_{n,n}$ and the
quotient space $\mathcal{P}_{n,n} / \mathcal{P}_{n,n-1}$ is one
dimensional generated by the n-degree wheel. If n is odd, then
$\mathcal{P}_n = \mathcal{P}_{n,n-1}$ (with the exception of
$\mathcal{P}_1$ being generated by $\eps{4mm}{theta}$)
\end{lemma}

\begin{proof}
See \cite{chm1,chm2,ba3} for proof and much more details on the
subject. Again, we mention that a wheel diagram looks like that: the
two wheel $\eps{7mm}{2wheel}$ , the 4-wheel $\eps{7mm}{4wheel}$ and
so on.\\

\qed
\end{proof}

This lemma transfers over to $\overrightarrow{\mathcal{A}(S^1)}$
where the primitives are the primitives of $\mathcal{A}(S^1)$
directed in all possible ways. Combining the previous two lemmas we
get:

\begin{corollary}
The only primitive diagrams that might possibly be non zero are the
wheel diagrams and $\eps{4mm}{theta}$.
\end{corollary}

However,

\begin{lemma}
The wheel diagrams are zero in $\overrightarrow{\mathcal{A}(S^1)}$.
\end{lemma}

\begin{proof}
We start by looking at a directed wheel diagram. There are only two
directed wheel diagrams of each degree, one for each of the two ways
of directing the inner loop, and it does not matter which one we
choose to look at. Due to the ``legal'' directing rules, it is not
hard to see that the external legs are always pointing towards the
base loop. We apply a directed STU relation on one of these legs
$\eps{7mm}{udu}=\eps{7mm}{sdu}-\eps{7mm}{tdu}$ obtaining two tree
diagrams (i.e. diagrams which contain no internal loops). A tree
diagram will always have exactly one leg pointing from the base loop
inwards. The tree diagram that corresponds to the first summand on
the right hand side of the relation has a vertex which looks like
$\eps{7mm}{zero}$. This means the tree diagram is equal to zero (see
the comment at the end of section 2.3). The second tree diagram,
corresponding to the second summand on the right hand side of the
directed STU relation, is ``leg crossed''. This crossing can be
untangled by moving the crossed leg around the base circle using the
directed AS relation when necessary and the commutativity of outward
pointing legs. Moving the leg around the circle will yield a diagram
which again has a vertex which looks like $\eps{7mm}{zero}$ (up to a
sign) and thus equals to zero. Apply this argument to any n-wheel
diagram and
the lemma is proven.\\
We demonstrate the entire argument on the 4-wheel :
\[
\eps{30mm}{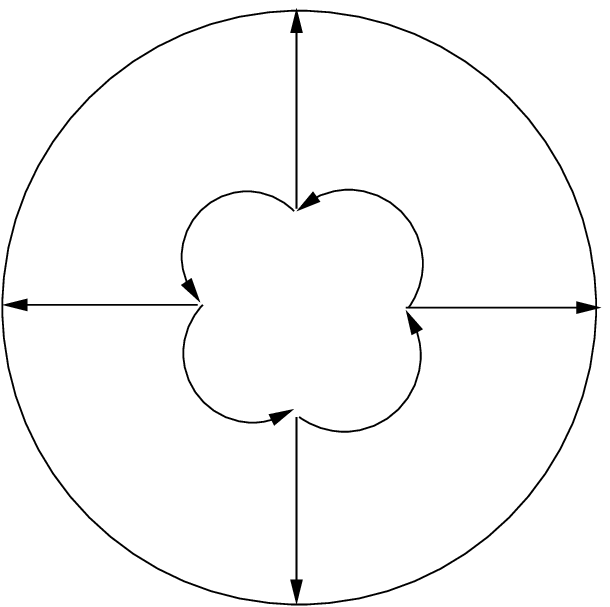} = \eps{30mm}{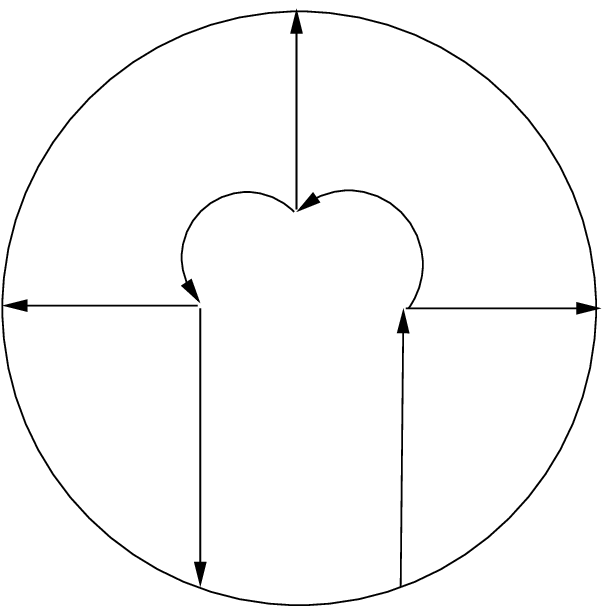} - \eps{30mm}{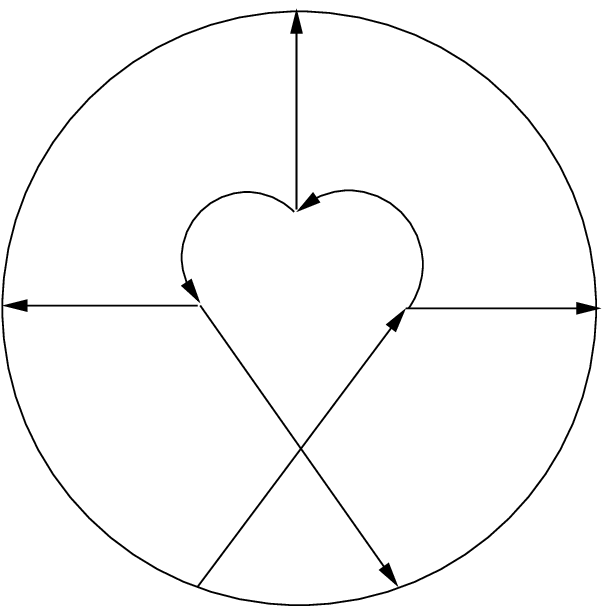} =0
-\eps{30mm}{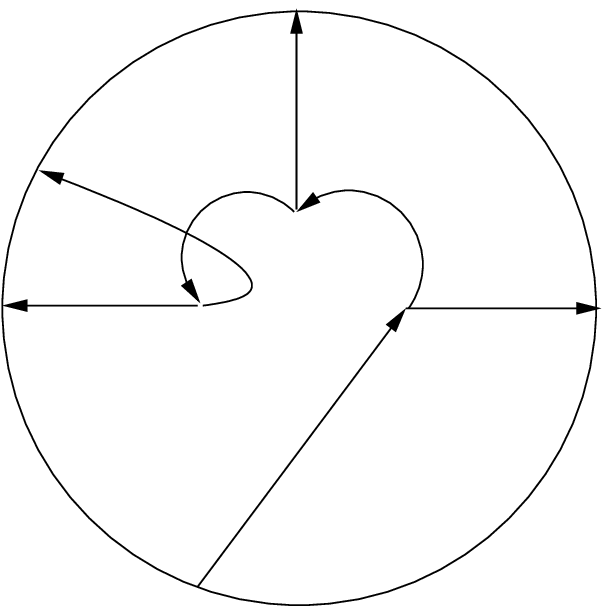} = 0 + 0 = 0
\]
\qed
\end{proof}

\vspace{3mm} Combining all of the above lemmas we finally get :

\begin{proposition}
The only possible non zero primitive element in
$\overrightarrow{\mathcal{A}(S^1)}$ is the directed single-chord
diagram $\eps{4mm}{theta}$. Thus, the only diagrams of degree $n$
that are possibly non zero in $\overrightarrow{\mathcal{A}(S^1)}$
are directed $\eps{4mm}{theta}^n$.
\end{proposition}

\begin{proof}
The proof is a corollary of all the above lemmas. The only primitive
element which is non-zero is $\eps{4mm}{theta}$ and it generates
powers of itself and sums of such powers.\qed
\end{proof}

\vspace{3mm} In CS perturbation theory it is well known that the
single-chord diagram has the framing (self-linking) of the knot as
its integral $\zeta^{CS}(\eps{4mm}{theta})$ (see for example
~\cite{gua1}). Higher degree diagrams, which are merely power of the
single-chord diagram, will hold information encoded in powers of the
framing of the knot. Moreover, this framing number can be normalized
as we wish, including normalization to zero, thus we can summarize
this chapter with the following proposition:

\begin{proposition}
CS topological quantum field theory with an inhomogeneous gauge
group and the standard knot observable (expectation value of the
trace, in some representation, of the holonomy along the knot) holds
no more than the framing information (which can be trivialized
anyway). \qed
\end{proposition}

CS theory with an inhomogeneous gauge group and the standard
observable fail to recognize knotting in $S^3$. As a gauge theory,
CS theory with inhomogeneous gauge group is equivalent to pure BF
theory and should therefore ``see'' knots at least for ISU(2) gauge
group, as was shown by Cattaneo et al. ~\cite{ca1,ca2,ca3}.
Therefore, we need a procedure that will extract some non trivial
information regarding the knot in our setting. This is what we will
be presenting in the next chapter.

\section{\bf{Extracting non trivial information regarding the knot}}

\subsection{Breaking gauge symmetry at one point - puncturing the space}
As seen in the previous chapter, the reason we got zero information
when the standard observable is used is the ``emptiness'' of
$\overrightarrow{\mathcal{A}(S^1)}$. Working with the inhomogeneous
gauge group allows us to factor through the ``legally'' directed
diagrams space when applying perturbation theory thus there is no
way of avoiding the consequences of lemmas 2.3 and 2.4 of the
previous chapter which state that the only possible non-vanishing
contribution in perturbation theory comes from the wheel diagrams.
The best we can do is try and avoid the reasons lemma 2.5 is true.
When we look at the proof of lemma 2.5 we can realize that the cause
for the vanishing contribution of the wheel diagrams is the fact
that legs can be moved around the Wilson loop (the base circle of
the diagram). Dealing with the loop translates (observable-wise) to
dealing with the trace (and its cyclic properties) in the standard
observable. In diagrammatic terms (perturbation-wise) this means
trying to cut the base circle open.
\\

Assume first that we just remove the trace from the standard
observable of a knot in $S^3$. We get an expectation value which is
formally a tensor in the universal enveloping algebra $\int{
\mathcal{DA} ~e^{iCS(A)}hol} \in U(L_0)$. This expression, though,
is not well defined since the holonomy by itself is not well defined
and transforms non-trivially by conjugation under gauge
transformation $\chi(x)$. Choosing (in advance) one point on the
knot, say $x_0$, the conjugation is done by the gauge element at
that point ($hol \rightarrow \chi(x_0)hol\chi^{-1}(x_0)$). The
remedy for that would be to view this (not well defined) invariant
in the co-invariant quotient $U(L_0)_{L_0}$. The way this is done
(observable-wise) is
of course by bringing back the trace. We have done nothing then.\\

This leads us to try and have special considerations for $x_0$ (the
chosen point on the knot) and the gauge symmetry at that point. Let
us break the gauge symmetry at $x_0$. There are a few equivalent
ways of looking at that process. The first is by declaring that the
gauge transformation at $x_0$ is always the identity element. This
might seem to be a bit artificial but an equivalent way of doing it
is taking $x_0$ out of the space on which the theory is defined by
puncturing $S^3$ on $x_0$. In other words, ``send'' $x_0$ to
infinity where gauge transformations are taken to vanish. This
shifts us to a theory on a long knot (embedding of $\mathbb{R}$) in
$\mathbb{R}^3$ with an invariant $\int{ \mathcal{DA} ~e^{iCS(A)}hol}
\in U(L_0)$ which is well defined in the framework of CS theory.\\

\emph{To summarize, we move from a CS theory on $S^3$ and a knot
(embedding of $S^1$) in it to a CS theory on $\mathbb{R}^3$ and a
long knot (embedding of $\mathbb{R}$) in it. This is done by
choosing a special point on the knot, puncturing the space there and
declaring the puncture as infinity. The standard observable
$Tr_Rhol$ is replaced by a more general observable $hol$.}\\

As far as knots are concerned, ``regular'' knots and long knots are
equivalent and we do not lose any information regarding the knot in
the process above. Moreover, $hol$ (and functions of it) is actually
the most general observable in this framework since the connection
can always be recovered (up to gauge transformations) from the
holonomy information. Thus we work in the setting which allows us to
extract maximum information regarding the knot when using the
inhomogeneous gauge group.
\\

We take now $\int{ \mathcal{DA} ~e^{iCS(A)}hol}$ and apply the same
perturbation theory factorization as in chapter 2. There are several
differences to note:
\\

 First, we work with a long knot. Thus, all the Feynman
diagrams are not based on a loop but on an interval (representing
the embedding of our long knot $\mathbb{R}$). The other properties
of the diagrams do not change, though. We define :

\begin{definition}
$\mathcal{A}(I)$ is the algebra of all connected trivalent graphs
based on an interval, modulo the STU, IHX and AS relations. (i.e.
$\mathcal{A}(S^1)$ with the circle replaced by an interval).\\
$\overrightarrow{\mathcal{A}(I)}$ is defined the same way as
$\overrightarrow{\mathcal{A}(S^1)}$, except that the diagrams are
based on an interval and not on a circle.
\end{definition}

Second, the knot invariant factors through the invariant subspace of
the universal enveloping algebra $U(L_0)^{L_0}$. This time there is
no necessity in taking the co-invariant quotient of that space, a
necessity that came about because of the base circle (or the trace
in the observable). Note that the fact that we indeed get a tensor
in the invariant subspace tells us that the construction was
independent of the choice of the special point $x_0$.\\

Third, we do not have a map to $\mathbb{C}$ yet.\\

All together we get the following factorization for the perturbation
theory of CS theory on $\mathbb{R}^3$ with a long knot in it (we
denote by $l$ the last map):

\begin{equation}
\{knots\} \rightarrow \mathcal{A}(I)[[k^{-1}]] \rightarrow
\overrightarrow{\mathcal{A}(I)}[[k^{-1}]] \xrightarrow{l}
  U(L_0)^{L_0}[[k^{-1}]]
\end{equation}
\\

\subsection{Building an observable and extracting non trivial information}
We have not shown so far that the algebra
$\overrightarrow{\mathcal{A}(I)}$ is non-zero and that our
construction is indeed non-trivial in the sense that it fixes the
problems encountered with the standard observable. We have also not
defined a scalar observable, i.e. a map $U(L_0)^{L_0} \rightarrow
\mathbb{C}$ which completes the perturbation theory factorization.
We will now achieve both goals by finding an explicit functional
from $\overrightarrow{\mathcal{A}(I)}$ to $\mathbb{C}$ that does not
vanish when composed with the above factorization. We will give a
formula for a scalar function $U(L_0)^{L_0} \rightarrow \mathbb{C}$
that is not zero on the part of $U(L_0)^{L_0}$ which comes from
wheel diagrams in $\overrightarrow{\mathcal{A}(I)}$. This observable
cannot possibly be any trace of $hol$ in some representation, since
that will bring us back to the case of chapter 2. Thus we need to
find some other scalar function of $hol$.\\

Start with an n-dimensional representation of $\mathcal{G}$ and
denote it $B$. Define the following representation $R$ of $L_0$:
\begin{equation}
R(g)= \biggl(
\begin{matrix} B(g) & 0
\\ 0 & B(g) \end{matrix} \biggl), \hspace{2mm} g \in \mathcal{G}
\hspace{20mm} R(\overline{g})= \biggl(
\begin{matrix} 0 & B(g) \\ 0 & 0 \end{matrix} \biggl), \hspace{2mm} \overline{g} \in \overline{\mathcal{G}}
\end{equation}
It is straight forward to verify that the commutation relation for R
is indeed compatible with the commutation relations of the
semi-direct product and that $R$ is indeed a representation of
$L_0$. We extend $R$ to $U(L_0)$ in the usual way.
\\

\begin{definition} ~\\

1. Let $\Delta$ denote the co-multiplication of the universal
enveloping algebra $U(L_0)$ and let $\Delta^m$ denote the map
$U(L_0) \rightarrow U(L_0)^{\otimes m}$ we get by composing $m-1$
times the map $\Delta$.
\\

2. Given the representation $R$ above, denote by $R^{\otimes m}$ the
m-tensored representation $U(L_0)^{\otimes m} \rightarrow
End(\mathbb{C}^{2n})^{\otimes m}$ defined by $u_1 \otimes \ldots
\otimes u_m \mapsto R(u_1) \otimes \ldots \otimes R(u_m)$.\\

3. Define a transformation $\lambda$: $End(\mathbb{C}^{2n})^{\otimes
m} \rightarrow End(\mathbb{C}^{2n})$ by $\tau_1 \otimes \tau_2
\otimes \ldots \otimes \tau_m \longmapsto \tau_1 \cdot C \cdot
\tau_2 \cdot C \ldots \cdot \tau_m \cdot C$, where $C$ is the
following $2n\times2n$ matrix: $\biggl(
\begin{matrix} 0 & 0 \\ I & 0
\end{matrix} \biggl)$ and $\cdot$ is matrix multiplication.\\

4. The ``half-trace'', denoted $Tr_{\frac{1}{2}}$, is the trace over
the upper-left $n\times n$ block of a $2n\times2n$ matrix:
\[
Tr_{\frac{1}{2}}\biggl(
\begin{matrix} A & B
\\ C & D \end{matrix} \biggl) = Tr(A)
\]
\\

\end{definition}

Armed with the above definitions and notations we can finally build:

\begin{definition}
Let $\Sigma_m$ denote the following composition of maps ($\circ$ is
used for composition):
\[
Tr_{\frac{1}{2}} \circ \lambda \circ R^{\otimes m} \circ \Delta^m :
U(L_0)^{L_0} \xrightarrow{\Sigma_m} \mathbb{C}
\]
\end{definition}

\begin{lemma}
The map $\Sigma_m \circ l$ is non-zero on the m-wheel diagram for at
least one choice of gauge group $L_0$ and representation $R$.
\end{lemma}

\begin{proof}
We start by looking at a directed wheel diagram (i.e. a wheel in
$\overrightarrow{\mathcal{A}(I)}$). All the external legs (the ones
that touch the base interval) point towards the interval, as we have
already seen. This means that after applying $l$, the tensor $w$
that we get in $U(L_0)$ will have all its components in
$\overline{\mathcal{G}}$. After applying $\Delta^m$ to $w$ we get a
tensor in $U(L_0)^{\otimes m}$ such that each of its summands
actually looks like $u_1 \otimes \cdots \otimes u_m \in
U(\overline{\mathcal{G}})^{\otimes m} \subset U(L_0)^{\otimes m}$.
\\

Now, the only type of tensors $u_1 \otimes \cdots \otimes u_m \in
U(\overline{\mathcal{G}})^{\otimes m}$ that will possibly not result
in zero after applying $\lambda \circ R^{\otimes m}$ to it, is the
type in which each component $u_i$ is actually an element of
$\overline{\mathcal{G}}$ (i.e. tensors of length one in $U(L_0)$
with entry from $\overline{\mathcal{G}}$). This is immediate from
the definition of $R$ (4) on $U(\overline{\mathcal{G}})$, the
definition of $\lambda$ and the definition of $C$.
\\

We summarize then -- given a wheel diagram we apply $l$ to it and
get a tensor $w$ in $U(L_0)$. $\Delta^m(w)$ is a tensor in
$U(L_0)^{\otimes m}$ but the only non zero contributions to $\lambda
\circ R^{\otimes m}(\Delta^m(w))$ come from the summands of
$\Delta^m(w)$ that look like  $w$ or any permutations of its entries
:

\[
w=l(directed ~m-wheel)= \sum_{i_1,\ldots
,i_m=1}^{dim\mathcal{G}}{C^{i_1,\ldots ,i_m}
\overline{X}_{i_1}\otimes\cdots \otimes \overline{X}_{i_m}}
\]

\[
\lambda \circ R^{\otimes m} (\Delta^m(w)) = \lambda \circ R^{\otimes
m} (\sum_{i_1,\ldots ,i_m=1}^{dim\mathcal{G}} \sum_{\sigma \in
S_m}{{C^{i_1,\ldots ,i_m}
\overline{X}_{i_{\sigma(1)}}\bigotimes\cdots \bigotimes
\overline{X}_{i_{\sigma(m)}}}})
\]

Where we made a distinction between the tensor product $\otimes$ of
$U(L_0)$ and the product $\bigotimes$ of $U(L_0)^{\otimes m}$ for
the sake of clarity.
\\

The calculation of $Tr_{\frac{1}{2}} \circ  \lambda \circ R^{\otimes
m} \circ \Delta^m \circ l(directed ~m-wheel)$ is almost the same as
the calculation of the weight system $W_{\mathcal{G},B}$ on the non
directed $m$-wheel diagram (i.e. the directed $m$-wheel with its
arrows forgotten). The only difference is the summation over all
permutations of the components of the tensor $w$ above, but that can
be done on the diagram level by permuting all external legs of the
non directed diagram. Thus we have :
\[
Tr_{\frac{1}{2}} \circ \lambda \circ R^{\otimes m} \circ \Delta^m
\circ l(directed ~m-wheel) = W_{\mathcal{G},B}(\chi(m-wheel))
\]
where $\chi(D)$ is defined to be the sum over all diagrams which
differ from $D$ by permutation of the external legs.\\

In order to show that $\Sigma_m \circ l(directed ~m-wheel)$ is
non-zero for at least one choice of $L_0$ and representation $R$ we
let $L_0=gl(n) \bigoplus \overline{gl(n)}$ with $B$ taken to be the
defining representation and calculate the highest order in $n$ of
$W_{gl(n),B}(\chi(m-wheel))$. Using simple counting arguments one
may show that all the highest order contributions are positive and
sum up to a non-zero contribution : $\Sigma_m \circ l(directed
~m-wheel) = mn^{m+1}$ + lower order terms (in $n$). We leave the
exact calculation and arguments to the interested reader. These can
be done using exercise 6.36 in \cite{ba3}.

\qed

\end{proof}

\begin{proposition}
In $\overrightarrow{\mathcal{A}(I)}$, the wheel diagrams (and
possibly the directed $\eps{4mm}{theta}$) are non zero.
\end{proposition}

\begin{proof}
We already know that the only possible non-zero directed generators
in $\overrightarrow{\mathcal{A}(I)}$ are the wheel diagrams and
possibly $\eps{4mm}{theta}$ (following the exact same arguments for
$\overrightarrow{\mathcal{A}(S^1)}$ as in chapter 2). Applying the
lemma above for all possible $m$ implies the proposition.

\qed

\end{proof}

\vspace{3mm} \textbf{The actual observable:} The proposition and
lemma above tell us that the perturbation theory (3) for the long
knot in $\mathbb{R}^3$ factors through a non-trivial space, and that
there is a way of extracting non-trivial information regarding the
knot. We can now build an actual quantum observable (function of
$hol$) to plug in the path integral such that the map
$Tr_{\frac{1}{2}} \lambda R^{\otimes m} \Delta^m$ will be the last
part of the factorization ($U(L_0)^{L_0} \rightarrow \mathbb{C}$) in
its perturbation theory calculation of the expectation value. This
observable is $Tr_{\frac{1}{2}}((R(hol) \cdot C)^m)$, where $R(hol)$
stands for representing the holonomy using the representation $R$ as
defined above.
\\

The expectation value $\int\mathcal{DA}~{e^{iCS}~Tr_{1/2}((R(hol)
\cdot C)^m)}$ has a perturbation theory factorization (with
$\Sigma_m$ as the last map):
\\
\begin{equation}
\{knots\} \rightarrow \mathcal{A}(I)[[k^{-1}]] \rightarrow
\overrightarrow{\mathcal{A}(I)}[[k^{-1}]] \xrightarrow{l}
U(L_0)^{L_0}[[k^{-1}]] \xrightarrow{\Sigma_m} \mathbb{C}[[k^{-1}]]
\end{equation}
\\

\subsection{The most general knot invariant in this setting - The Alexander-Conway polynomial}
We have just seen on the last section that we can extract all the
information contained in the wheel diagrams by using CS theory with
the inhomogeneous gauge group on long knots in $\mathbb{R}^3$ (with
the new class of observables). This is opposed to the fact that the
standard observable in CS theory with the inhomogeneous gauge group
will not see these knots in $S^3$. We now turn to the task of
recognizing the non trivial information we extracted.

\begin{theorem}
The most general knot invariant coming from Chern-Simons topological
quantum field theory with an inhomogeneous gauge group is
non-trivial and is in the algebra generated by the coefficients of
the Alexander-Conway polynomial (together with possible framing
information).
\end{theorem}

\begin{proof} In \cite{kri2}, following motivations from \cite{kri1,ba4}, it is
proven that any weight system which is zero on the kernel of the
de-framing operator and zero on all the primitive spaces which are
not wheels, belongs to a knot invariant which is in the algebra
generated by the coefficients of the Alexander-Conway polynomial.
Since we showed that the only diagrams that contribute to our weight
system are the wheel diagrams, and since we can always choose our
knot to be of (standard) zero framing, the theorem follows. Other
choices of framing will show up through $\eps{4mm}{theta}$.

\qed

\end{proof}

\vspace{3mm} Notice that it is possible to take the semi-direct
product of $G$ with $\mathcal{G}^*$, the dual of the Lie algebra. We
can follow the exact same considerations and get the same results,
this time without the need of the metric on $\mathcal{G}$ (just use
the natural dual pairing to get a metric on $L_0$). Therefore our
results are true for any Lie algebra $\mathcal{G}$, not necessarily
a metrizable one.

\section{\bf{Pure BF Theory and the Alexander-Conway polynomial}}
Recall that the pure BF theory has the following action $S_{BF} =
\int_{S^3} Tr(B \wedge F)$, where F is the curvature of the
G-connection A and B is a one form taking values in the algebra
$\mathcal{G}$. In \cite{ca1,ca2,ca3} Cattaneo et al. took
$Tr_R(Exp(\lambda\gamma_1))$ as the most general observable one can
consider in Pure BF theory, as long as the standard (zero) framing
of the knot is chosen. Here R is any representation of G, $\gamma_1$
is $\oint_x{Hol(A)B(x)Hol(A)}$ (the first degree element in the
Taylor expansion of $Hol(A + \kappa B)$ around $\kappa =0$) and
$\lambda$ is a ``coupling'' constant counting orders.\\

Using this observable, QFT techniques and the Melvin-Morton
conjecture (MMR, see \cite{ba4}), it was proven \cite{ca1} that the
set of unframed knot invariants we get from the pure BF theory,
using SU(2) gauge group, is generated by the coefficients of the
Alexander-Conway polynomial.
\\

In this paper we proved:

\begin{theorem}
The most general knot invariant coming from pure BF topological
quantum field theory with any gauge group whose Lie algebra is
metrizable and with any representation, is in the algebra generated
by the coefficients of the Alexander-Conway polynomial.
\end{theorem}

\begin{proof}
Using the fact that the pure BF theory with gauge group G is just CS
theory with gauge group IG ~\cite{bir2} this theorem is just a
corollary of Theorem 1.

\qed
\end{proof}

\vspace{3mm} The same generalization for non-metrizable Lie
algebras, as right after Theorem 1, applies here as well.\\

A closer look at our observable $Tr_{\frac{1}{2}}((R(hol)\cdot
C)^m)$ shows that it is actually the $m^{th}$ degree in the
expansion of the BF observable $Tr_R(Exp(\lambda\gamma_1))$ (up to
numerical factors coming from the exponent expansion). The fact that
the BF observable also breaks the gauge symmetry of the pure BF
theory at one point (it has to be assumed that the special BF gauge
symmetry is identity at one point on the knot) was not emphasized
much before (though observed of course in \cite{ca1,ca3,ca4}). Our
setting gives a natural explanation as to why this is indeed the
most general observable for BF theory without any need for taking
various limits or referring to MMR. We also get a somewhat clearer
explanation as to why one can ignore the $\gamma_0$ part of the
Taylor expansion above.
\\

One can now say that as far as knot theory and 3 dimensional BF
topological quantum field theory (with or without cosmological
constant) are concerned, there is nothing new beyond Chern-Simons
theory, which can reproduce the same knot theoretical results.

\section{\bf{Discussion: Non-Semi-Simplicity}}
Algebraically, the main difference between our construction and the
standard one is getting invariants in the invariant subspace
$U(L_0)^{L_0}$ instead of the co-invariant quotient space
$(U(L_0)^{L_0})_{L_0}$ of the universal enveloping algebra.
\\

Our construction gives no new extra information if we work with a
semi-simple gauge group instead of $L_0$. This is true because for
semi-simple groups, $U(\mathcal{G})^{\mathcal{G}}$ is isomorphic to
$(U(\mathcal{G})^{\mathcal{G}})_{\mathcal{G}} \cong
U(\mathcal{G})_{\mathcal{G}}$. We get no new information from
applying our procedure and working in $\mathcal{A}(I)$ instead of
$\mathcal{A}(S^1)$, since these spaces are isomorphic (the
diagram-space way of expressing the previous isomorphism). Our
construction seems to detect the difference between the invariant
subspace and the co-invariant quotient of the non-semi-simple Lie
algebra we have worked with. In other words, it uses the fact that
$\overrightarrow{\mathcal{A}(I)}$ is not isomorphic to
$\overrightarrow{\mathcal{A}(S^1)}$ (the map
$\overrightarrow{\mathcal{A}(I)} \rightarrow
\overrightarrow{\mathcal{A}(S^1)}$ that closes the base interval
into the base circle has a kernel!).
\\

We saw that there is more to observe in Chern-Simons theory when
dealing with a non semi-simple gauge group. A question to be raised
is whether other non semi-simple gauge groups can hold the same
property as the inhomogeneous group - meaning that the above
procedure can extract more information about the knot than what was
given by the standard observables.

\section{\bf{Acknowledgments}}
I wish to thank Dror Bar-Natan for his endless support in this
research. I also wish to thank Alberto Cattaneo for reading an early
version of this paper and giving his comments.

\end{document}